\documentclass[a4paper,12pt]{article} 
\newdimen\paperhight
\usepackage{amsmath}
\usepackage{amssymb}
\usepackage{amsfonts}
\usepackage[usenames]{color}

\newcommand{\wt}{{\rm wt}}

\newcommand{\pr}{\par \vspace{3mm}\noindent [{\bf Proof}] \qquad}
\newcommand{\prend}{\hfill \qed \par \vspace{3mm}}

\newcommand{\qed}{\quad\hbox{\rule[-2pt]{3pt}{6pt}}\par\vspace{3mm}}

\newcommand{\C}{\mathbb C} 
\newcommand{\Z}{\mathbb Z}

\newcommand{\N}{\mathbb N}

\newcommand{\CI}{{\cal I}}
\newcommand{\CJ}{{\cal J}}
\newcommand{\CF}{{\cal F}}

\newcommand{\CO}{{\cal O}}

\newcommand{\CY}{{\cal Y}}

\newcommand{\Hom}{{\rm Hom}}

\newcommand{\Ker}{{\rm Ker}}

\begin{document}
\title{$C_1$-cofiniteness and Fusion Products of Vertex Operator Algebras}
\author{\begin{tabular}{c}
Masahiko Miyamoto \\
Institute of Mathematics, \\
University of Tsukuba, \\
Tsukuba, 305 Japan \end{tabular}}
\date{}
\maketitle

\begin{abstract}
Let $V$ be a vertex operator algebra. 
We prove that if $U$ and $W$ are $C_1$-cofinite $\N$-gradable $V$-modules, then 
a fusion product $U\boxtimes W$ is also a $C_1$-cofinite $\N$-gradable $V$-module, 
where the fusion product is defined by (logarithmic) intertwining operators. 
\end{abstract}

\section{Introduction}
The tensor product theory is a powerful tool in the theory of 
representations.  
Unfortunately, in the theory of vertex operator algebras (shortly VOA), 
a tensor product (we call ``fusion product'') for some modules 
may not exist in the category of modules of 
vertex operator algebras. In order to avoid such an ambiguity, 
we will introduce a new approach to treat fusion products. 
Let us explain it briefly. The details are given in \S 3. 
Let $V=\oplus_{n=K}^{\infty}V_n$ be a vertex operator algebra 
(shortly VOA) and ${\rm mod}_{\N}(V)$ denote 
the set of $\N$-gradable (weak) $V$-modules, 
where a (weak) $V$-module $W$ is called $\N$-gradable if 
$W=\oplus_{m=0}^{\infty} W_{(m)}$ such that 
$$\displaystyle{v_kw\in W_{(m+\wt(v)-k-1)}}$$ 
for any homogeneous element $v\in V_{\wt(v)}$, $k\in \Z$ and $w\in W_{(m)}$.  
It is well-known that 
$$g(V):=V\otimes_{\C}\C[x,x^{-1}]/(L(-1)\otimes 1-1\otimes\frac{d}{dx})V\otimes_{\C}\C[x,x^{-1}]$$ 
has a Lie algebra structure and all (weak) $V$-modules are $g(V)$-modules (see \cite{B}). 
For $U,W\in {\rm mod}_{\N}(V)$, we introduce 
a $g(V)$-module $U\boxtimes W$ (or its isomorphism class) 
as a projective limit of a direct set of $V$-modules (by viewing them as $g(V)$-modules). 
So, a $g(V)$-module $U\boxtimes W$ always exists. The key point is 
that a fusion product for $U,W\in {\rm mod}_{\N}(V)$ exists if and only if 
$U\boxtimes W$ is a $V$-module (and so it is a fusion product). 

The main purpose of this paper is to explain the 
fusion products by emphasizing the importance of $C_1$-cofiniteness. 
The importance of the $C_1$-cofiniteness conditions on modules was firstly 
noticed by Huang in \cite{H}, 
where he has proved that intertwining operators of a $C_1$-cofinite 
$\N$-gradable module 
from a $C_1$-cofinite $\N$-gradable module to an $\N$-gradable 
module satisfy  
a differential equation. He has also shown the associativity 
of intertwining operators among $C_1$-cofinite $\N$-gradable  
modules by using the space of solutions of this differential equation.  
We will prove Key Theorem by using his idea. 
In order to follow his arguments, we give a slightly 
different definition of $C_m$-cofiniteness for modules. \\

\noindent
{\bf Definition 1} \quad {\it 
Set $m=1,2,\ldots$. 
A $V$-module $U$ is said to be ``$C_m$-cofinite as a $V$-module'' 
if $C_m(U):= {\rm Span}_{\C}\{v_{-m}u \mid u\in U, v\in V, \wt(v)>1-m\}$ 
has a finite codimension in $U$.}\\

This is slightly different from the old one. 
For example, any VOA $V$ is always $C_1$-cofinite as a $V$-module 
in our definition. 
Since $(L(-1)v)_{-m}=mv_{-m-1}$ and $\wt(L(-1)v)=\wt(v)+1$, 
$C_m$-cofiniteness implies $C_{m-1}$-cofiniteness for $m=2,3,\ldots.$

We will prove the following theorem. \\

\noindent
{\bf Key Theorem} \quad {\it
Let $V$ be a VOA. 
For each $m=1,2,\ldots$ and $C_m$-cofinite $\N$-gradable $V$-modules $U$ and $W$, 
there is an integer $f_m(U,W)$ such that 
if $T$ is an $\N$-gradable $V$-module and 
there is a surjective (logarithmic) intertwining operator in $I\binom{T}{U\, \, W}$,  
then $\dim T/C_m(T)<f_m(U,W)$. 
In particular, $T$ is also $C_m$-cofinite as a $V$-module. }\\

Hereafter $I\binom{T}{U\, \, W}$ denotes the space of (logarithmic) intertwining operators of 
type $\binom{T}{U\, \, W}$. 
For $\phi\in {\rm Hom}(T,S)$ and $\CY\in I\binom{T}{U\,\,W}$, we can define an 
intertwining operator $\phi\circ\CY$ of type $\binom{S}{U\,\,W}$ by 
$$(\phi\circ\CY)(u,z)w:=\phi(\CY(u,z)w) \quad \mbox{ for }u\in U, w\in W.$$ 
$\CY\in I\binom{T}{U\, \, W}$ is called 
``surjective'' if for any proper injection $\epsilon:E\to T$ and 
any intertwining operator $\CJ\in I\binom{E}{U\, \, W}$, 
we have $\epsilon\circ\CJ\not=\CY$. 
As an application, we will prove the following theorem. \\

\noindent
{\bf Main Theorem} \quad {\it
Let $m=1,2,\ldots$ and let $V$ be a VOA.  
If $U$ and $W$ are $\N$-gradable $C_m$-cofinite $V$-modules, then 
a fusion product $U\boxtimes W$ is also a $C_m$-cofinite $\N$-gradable $V$-module. }\\

As an application of Main Theorem, we have: \\

\noindent
{\bf Corollary 2} \quad 
{\it Let $V$ be a simple VOA with $V\cong V'$, where $V'$ denotes 
the restricted dual of $V$.  
If there is a $V$-module $W$ such that 
$W$ and its restricted dual $W'$ are both $C_2$-cofinite, 
then $V$ is $C_2$-cofinite. } \\

\section{Proof of Key Theorem}
We first assume that $U$ and $W$ are indecomposable. 
Since $U$ and $W$ are $C_m$-cofinite, there is $N\in \N$ such that  
$U=C_m(U)+E$ and $W=C_m(W)+F$, where 
$E=\oplus_{k=0}^NU_{r_U+k}$ and $F=\oplus_{k=0}^NW_{r_W+k}$ 
and $r_U$ and $r_W$ denote the lowest weights of 
$V$-modules $U$ and $W$, respectively.  
We fix bases $\{p^i\mid i\in I \}$ of $E$ 
and $\{q^j\mid j\in J \}$ of $F$ consisting of homogeneous elements, 
respectively.   

Let $\CY\in I\binom{T}{U\quad W}$ be a surjective intertwining operator and 
let $T'$ denote the restricted dual of $T$. 
For each $\theta\in T', u\in U, w\in W$, we define  
a bilinear form 
$$\langle \theta, \CY(u,z)w\rangle\in \C\{z\}[\log z]$$ 
by $\theta(\CY(u,z)w)$. 
Applying the idea in \cite{H} to $\theta\in {\rm Annih}(C_m(T))\cap T'$, 
we have the following lemma. \\

\noindent
{\bf Lemma 3} \quad {\it
For $p\in U$, $q\in W$ and $\theta\in {\rm Annih}(C_m(T))\cap T'$, 
$$ F(\theta,p,q;z):=\langle \theta, \CY(p,z)q\rangle \vspace{-2mm}$$ 
is a linear combination of $\{F(\theta,p^i,q^j;z)\mid i\in I, j\in J\}$ with 
coefficients in $\C[z,z^{-1}]$. 
We are able to choose these coefficients independently from 
the choice of $\theta$. Moreover, there is an integer $f_m(U,W)$ given by $U$ and $W$ only  
such that $\dim(C/C_m(T)) <f_m(U,W)$.} \\

\pr
We will prove Lemma 3 for $m=1$. For $m\geq 2$, the proofs are similar. 
We will prove the first assertion in Lemma 3 by induction on $\wt(p)+\wt(q)$. 
Clearly, we may assume 
that $\wt(p)>N+r_U$ or $\wt(q)>N+r_W$. 
If $\wt(p)>N+r_U$, then $p=\sum_k v^k_{-1}a^k$ 
with $v^k\in V$ and $a^k\in U$. 
We note that this expression does not depend on the choice of $\theta$. 
Since $p$ is a linear sum, we are able to treat each term separately, that is, 
we may assume $p=v_{-1}a$  with $v\in V$ and $a\in U$.  
Then for $\theta\in {\rm Annih}(C_1(T))\cap T'$, we have:
$$\begin{array}{rl}
\displaystyle{\langle \theta,\CY(p,z)q\rangle=}&\displaystyle{\langle\theta,\CY(v_{-1}a,z)q\rangle}\vspace{2mm}\cr
=&\displaystyle{\langle\theta, Y^-(v,z)\CY(a,z)q+\CY(a,z)Y^+(v,z)q\rangle }\vspace{2mm}\cr
=&\displaystyle{\langle\theta, \CY(a,z)Y^+(v,z)q\rangle},
\end{array}\eqno{(2.1)}$$
where $\displaystyle{Y^-(v,z)=\sum_{h<0}v_hz^{-h-1}}$ and 
$\displaystyle{Y^+(v,z)=\sum_{h\geq 0}v_hz^{-h-1}}$. 
This is a reduction on the sum of weights 
because $\wt(v_{h}q)< \wt(v)+\wt(q)$ for $h\geq 0$, that is, 
all terms of $Y^+(v,z)q$ have less weights than $\wt(v)+\wt(q)$. 
An important thing is that the processes of these reductions 
are irrelevant with the choice of $\theta$.

Similarly, if $\wt(q)>N+r_W$, then we may assume 
$q=v_{-1}b$ with $v\in V$ and $b\in W$ and 
we have:
$$\begin{array}{rl}
\langle\theta,\CY(p,z)q\rangle=&\langle \theta, \CY(p,z)v_{-1}b\rangle\vspace{2mm}\cr
=&\displaystyle{\langle\theta,v_{-1}\CY(p,z)b\rangle
+\langle \theta,\sum_{i=0}^{\infty}\binom{-1}{i}z^{-1-i} \CY(v_ip,z)b \rangle}\vspace{2mm}\cr
=&\displaystyle{\sum_{i=0}^{\infty}\binom{-1}{i}z^{-1-i}
\langle \theta, \CY(v_ip,z)b \rangle.}
\end{array}\eqno{(2,2)}$$
We also note that this calculation is independent of the choice of 
$\theta$ and 
this is also a reduction on the weights because 
$\wt(v_ip)+\wt(b)<\wt(v_{-1}b)+\wt(p)$ for $i\geq 0$. 
Therefore, $\langle \theta, \CY(p,z)q\rangle$ 
is a linear combination of 
$\{\langle \theta, \CY(p^i,z)q^j\rangle \mid i\in I, j\in J\}$ 
with coefficients in $\C[z,z^{-1}]$ and these coefficients are 
independent of the choice of $\theta$. 

We next prove the second assertion in Lemma 3. 
We consider an $|I|\times |J|$-dimensional  
vector $A=(\langle \theta, \CY(p^i,z)q^j\rangle)_{i,j}$. 
Then there is a square matrix 
$B\in M_{|I|\times |J|,|I|\times |J|}(\C[z,z^{-1}])$ which 
does not depend on the choice of $\theta$ such that 
$$\frac{d}{dz}A=(\langle \theta, \CY(L(-1)p^i,z)q^j\rangle)_{i,j}=BA. \eqno{(2.3)}$$ 
The space of solutions of the above differential equation (2.3) 
has a finite dimension and so the space of the choice of $\theta$ 
in $A$ is also of finite dimension.   
Since $\CY$ is surjective, $\dim T/C_1(T)$ 
is bounded by a number $f_1(U,W)$ which does depend on $U$ and $W$ only. 
In particular, $T$ is $C_1$-cofinite.  

We next assume that $U=\oplus U^{(r)}$ and $W=\oplus W^{(k)}$ with 
indecomposable $V$-modules $U^{(r)}$ and $W^{(k)}$. 
Clearly, $U^{(r)}$ and $W^{(k)}$ are also $C_1$-cofinite as $V$-modules. 
Since $\wt(v_{-1}u)>\wt(u)$ for $v\in V$, $u\in U^{(r)}$ with $\wt(v)>0$, 
we have $U^{(r)}/C_1(U^{(r)})\not=0$ for every $r$ and so 
$U$ is a finite direct sum of indecomposable modules. Similarly, so is $W$.  
For $\CY\in I\binom{T}{U\, \, W}$ and for each $(r,k)$, we define 
$\CY^{(r,k)}\in I\binom{T^{(r,k)}}{U^{(r)}\, \, W^{(k)}}$ by restrictions, that is, 
$$\CY^{(r,k)}(u^{(r)},z)w^{(k)}:=\CY(u^{(r)},z)w^{(k)}\quad \mbox{ for }u^{(r)}\in U^{(r)}, 
w^{(k)}\in W^{(k)}  $$
and $T^{(r,k)}$ is the subspace of $T$ spanned by all coefficients of 
$\CY^{(r,k)}(u^{(r)},z)w^{(k)}$. 
As we showed, there are integers $f_1(U^{(r)},W^{(k)})$ such that 
$\dim T^{(r,k)}/C_1(T^{(r,k)})\leq f_1(U^{(r)},W^{(k)})$. 
Since $T=\sum_{r,k} T^{(r,k)}$ and $C_1(T^{(r,k)})\subseteq C_1(T)$, we have 
$$\dim T/C_1(T)\leq \sum_{r,k} \dim T^{(r,k)}/C_1(T^{(r,k)}) \leq \sum_{r,k} f_1(U^{(r)},W^{(k)})$$ 
as we desired. 

This completes the proof of Key Theorem. 
\prend

\section{On fusion products}
In this section, we would like to explain our approach to the fusion product of 
two modules.  
The fusion product of modules in the theory of vertex operator algebra 
are firstly defined by Huang and Lepowsky in several ways, (see \cite{HL}).  
We go back to the original concept of 
tensor products, that is, as stated in the introduction of \cite{HL}, 
it should be a universal one in the following sense, that is, 
if $U$ and $W$ are $V$-modules, then 
a fusion product is a pair $(U\boxtimes W, \CY^{U\boxtimes W})$ 
of a $V$-module $U\boxtimes W$ and 
an intertwining operator $\CY^{U\boxtimes W}\in I\binom{U\boxtimes W}{U\, \,W}$ 
such that for any $V$-module $T$ and any intertwining operator 
$\CY\in I\binom{T}{U\, \, W}$, there is 
a homomorphism $\phi:U\boxtimes W\to T$ such that 
$\phi\circ \CY^{U\boxtimes W}=\CY$, that is, 
$$\phi(\CY^{U\boxtimes W}(u,z)w)=\CY(u,z)w$$ 
for any $u\in U$ and $w\in W$. 
Unfortunately, in the theory of vertex operator algebra, 
unlike the categories of vector spaces, 
a fusion product module may not exist. 

Our idea in this paper is that we first construct a 
$g(V)$-module $U\boxtimes W$.  
Furthermore, as we will see, $U\boxtimes W$ satisfies all conditions (Commutativity, etc) 
to be a $V$-module except lower truncation property. 
Associativity is also true if it is well-defined, that is, 
if the lower truncation property holds on $U\boxtimes W$.  

Let us construct it. We fix $U,W\in {\rm mod}_{\N}(V)$ and  
consider the set of surjective intertwining operators $\CY$ of $U$ from $W$: 
$$ \CF(U,W)=\{(F,\CY)\mid  F\in {\rm mod}_{\N}(V), \CY\in I\binom{F}{U\, \, W}
\mbox{ is surjective}\}.$$
Here the set of intertwining operators includes not only 
formal $\C$-power series but also all intertwining operators of logarithmic forms. 

We define $(F^1,\CY_1)\cong (F^2,\CY_2)$ if there is an 
isomorphism $f:F^1\to F^2$ such that 
$f\circ\CY_1=\CY_2$, that is, $f(\CY_1(u,z)w)=\CY_2(u,z)w$ for any $u\in U, w\in W$. 
We also define a partial order $\leq $ in $\CF(U,W)/\cong$ as follows: \\
For $(F^1,\CY_1),(F^2,\CY_2)\in \CF(W,U)$, 
$$\CY_1\leq \CY_2 \Leftrightarrow {}^{\exists}f\in \Hom_V(F^2,F^1)  
\mbox{ such that }f\circ \CY_2=\CY_1.$$ 
We note that since $\CY_1$ and $\CY_2$ are surjective, $f$ is uniquely determined. 
Clearly, if $\CY_1\leq \CY_2$ and $\CY_2\leq \CY_1$, then we 
have $(F^1,\CY_1)\cong (F^2,\CY_2)$. \\

\noindent
{\bf Lemma 4} \quad {\it 
$\CF(U,W)/\cong$ is a (right) directed set. }\\

\pr 
For $\CY_1\in I\binom{F^1}{U\, \, W}$ and $\CY_2\in \CI\binom{F^2}{U\, \, W}$, 
we define $\CY$ by 
$$\CY(u,z)w=(\CY_1(u,z)w, \CY_2(u,z)w)\in (F^1\times F^2)\{z\}[\log z] \quad 
\mbox{ for }u\in U, w\in W.$$
Clearly, $\CY\in I\binom{F^1\times F^2}{U\, \, W}$. 
Let $F\subseteq F^1\times F^2$ be the subspace 
spanned by all coefficients of $\CY(u,z)w$ with $u\in U, w\in W$, 
then $(F,\CY)\in \CF(U,W)$. 
Moreover, by the projections $\pi_i:F^1\times F^2 \rightarrow F^i$, we have 
$\pi_1\circ\CY=\CY_1$ and $\pi_2\circ\CY=\CY_2$, that is, we have 
$(F^1,\CY_1)\leq (F,\CY)$ and $(F^2,\CY_2)\leq (F,\CY)$ as we desired. 
\prend

Since $\CF(U,W)/\cong$ is a direct set, we can consider 
a projective limit of $\CF(U,W)/\cong$ and we denote it (or 
a representative of its isomorphism class) by $(U\boxtimes W, \CY^{U\boxtimes W})$ 
(or simply $U\boxtimes W$). 
Since $T$ is a $g(V)$-module for any $(T,\CY)\in \CF(U,W)$, 
a projective limit $U\boxtimes W$ is also a $g(V)$-module. 
In order to see the actions of $g(V)$ on $U\boxtimes W$, 
let us show a direct construction of projective limit. 
Let $\{(F^i,\CY_i) \mid i\in I\}$ be the set of all representatives of $\CF(U,W)/\cong$ 
and set $\CY_i(u,z)=\sum_{j=0}^K\sum_{m\in \C} u^{\CY_i}_{(j,m)}z^{-m-1}\log^jz$ 
for $u\in U$. 
We consider the product $(\prod_{i\in I}F^i, \prod_{i\in I}\CY_i)$ and 
take subspaces $F_{r}$ of $\prod_{i\in I}F^i$ spanned by all coefficients 
$\prod_{i\in I} u_{j, {\rm wt}(w)-1-r+\wt(u)}^{\CY_i}w$ of weights $r\in \C$ 
for homogeneous elements $u\in U, w\in W$ and $j\in \N$. 
We then set 
$$F=\coprod_{r\in \C}F_{r} \qquad \left( \subseteq \prod_{i\in I}F^i\right). \eqno{(3.1)}$$ 
We note that $F$ is $\C$-gradable by the definition.

For $v^1,v^2\in V$, since the commutator formula  
$$[v^1_n,v^2_m]=\sum_{i=0}^{\infty}\binom{n}{i}(v^1_iv^2)_{n+m-i} \eqno{(3.2)}$$
holds on every $F^i$, we have the commutatior formula (3.2) on $\prod_{i\in I} F^i$ and also on $F$.  
$L(-1)$-derivative property is also true on $F$ since it is true on every $F^i$. 
For Associativity, since  
$$ (v^1_nv^2)_m=\sum_{i=0}^{\infty}\binom{n}{i}(-1)^i 
\{v^1_{n-i}v^2_{m+i}-(-1)^nv^2_{n+m-i}v^1_i\} \eqno{(3.3)}$$
is true for each $F^i$, it is also true on $F$ if the right side of (3.3) is well-defined on $F$. 
Therefore, if $V$ acts on $F$ truncationally, then $F$ becomes a (weak) $V$-module. 
 
Let us show $F\cong U\boxtimes W$ as $g(V)$-modules. 
Let $\pi_i:F\subseteq \prod_{h\in I}F^h \to F^i$ be a projection. By the definition, 
$\pi_i\circ(\prod_{h\in I}\CY_h)=\CY_i$. 
We note that if $(F^j,\CY_j)\leq (F^i,\CY_i)$ for $i,j\in I$, that is, if 
there is a surjective homomorphism $\phi_{i,j}:F^i \to F^j$ such that 
$\phi_{i,j}\circ\CY_i=\CY_j$, then 
$\pi_j\circ(\prod_{h\in I}\CY_h)=\CY_j=\phi_{i,j}\circ(\CY_i)
=\phi_{i,j}\circ(\pi_i\circ(\prod_{h\in I}\CY_h))$. 
Therefore, for any $\alpha\in F$, $\pi_j(\alpha)=\phi_{i,j}(\pi_i(\alpha))$ and 
so we have the following commutative diagram: 
$$ \begin{array}{rccc} 
\pi_i:&(F,\prod_{h\in I}\CY_h) & \longrightarrow & (F^i,\CY_i) \vspace{2mm}\cr
      &  || \mbox{ Identity}             &        & \downarrow \phi_{i,j} \vspace{2mm}\cr
\pi_j:&(F,\prod_{h\in I}\CY_h) &\longrightarrow & (F^j,\CY_j).\cr  
\end{array}$$
We note that since $\CY_i$ are surjective, $\phi_{i,j}$ are uniquely determined and 
$\phi_{jk}(\phi_{ij}(\CY_i(u,z)w))=\CY_k(u,z)w=\phi_{ik}(\CY_i(u,z)w)$.
Therefore, there is a surjective homomorphism $\phi:F\to U\boxtimes W$ such that 
$\phi\circ(\prod_{h\in I}\CY_h)=\CY^{U\boxtimes W}$. 
On the other hand, since $U\boxtimes W$ is a projective limit of $(F^i,\CY_i)_{i\in I}$, 
there are $\varphi_i:U\boxtimes W \to F^i$ such that $\varphi_i\circ(\CY^{U\boxtimes W})=\CY_i$. 
Therefore, $\varphi_i\circ(\phi\circ(\prod_{h\in I}\CY_h))=\CY_i=\pi_i\circ(\prod_{h\in I}\CY_h)$ and 
so $\varphi_i\phi=\pi_i$.  
Since $\cap_{h\in I}\Ker \pi_h=0$ by the definition, $\Ker\phi=0$ and so 
$(F, \prod_{h\in I}\CY^h)$ is isomorphic to a projective limit $(U\boxtimes W,\CY^{U\boxtimes W})$ 
of $\CF(U,W)/\cong $. \\

\noindent
{\bf Definition 5}\quad{\it 
We will call $U\boxtimes W$ a ``fusion product" of $U$ and $W$ (even 
if it is not a $V$-module.) }\\

We have to note that since $\CY^{U\boxtimes W}$ is a projective limit, 
the powers of $\log z$ in $\CY^{U\boxtimes W}(w,z)u$ may not have an upper bound 
even if $U\boxtimes W$ is a (weak) $V$-module. However, since $\CY^{U\boxtimes W}$ satisfies 
$L(-1)$-derivation, Commutativity and Associativity and etc, 
we still treat it as an intertwining operator.

\section{Proof of Main Theorem}
Before we start the proof of Main Theorem, 
we prove the following lemma. \\

\noindent
{\bf Lemma 6}\quad{\it 
If an $\N$-gradable module $U=\oplus_{m=0}^{\infty}U_{(m)}$ is 
$C_1$-cofinite as a $V$-module, then $\dim U_{(m)}< \infty$ 
for any $m\in \{0,1,\ldots\}$. }\\

\pr  We will prove it by induction on $m$. For $m=0$, 
since $U_{(0)}\cap C_1(U)=\{0\}$ by the definition of $C_1$-cofiniteness, 
we have $\dim U_{(0)}\leq \dim (U/C_1(U))<\infty$. 
We next assume $\dim U_{(k)}<\infty$ for $k=0,\ldots,m-1$. 
Then 
$$(V_i)_{-1}U_{(m-i)}:={\rm Span}_{\C}\{ v_{-1}u 
\mid v\in V_i, u\in U_{(m-i)}\}$$
is also of finite dimension for $i=1,2,\ldots,m$ since $\dim V_i<\infty$. 
Furthermore, since $U_{(m)}/\sum_{i=1}^m (V_i)_{-1}U_{(m-i)}\subseteq U/C_1(U)$, 
we have $\dim U_{(m)}<\infty$ as we desired.  
\prend

Let $U$ be a $C_1$-cofinite $\N$-gradable module. 
Since $\dim U_{(m)}<\infty$ and $L(0)$ acts on $U_{(m)}$, 
$U_{(m)}$ is a finite direct sum of generalized eigenspaces of $L(0)$. 
Furthermore, since $U$ is $C_1$-cofinite as a $V$-module, 
there is an integer $r$ such that 
$U_{(m)}\subseteq C_1(U)$ for $m\geq r$. 
Hence there is a finite set $\{a_i\in \C \mid i\in I\}$ such that  
$$  U=\oplus_{i\in I}\oplus_{m=0}^{\infty} U_{a_i+m}, $$
where $U_{a_i+m}$ is a generalized eigenspace of $U$ for $L(0)$ with eigenvalue $a_i+m$. 
We note $a_i\not\equiv a_j \pmod{\Z}$ for $i\not=j$. 
In particular, we have 
$$\wt(U)\subseteq \cup_{i\in I} (a_i+\N),$$
where $\wt(U)$ denotes the set of all weights of elements in $U$.  
We may assume $U_{(r)}=\oplus_{i\in I}U_{a_i+r}$ by rearranging the $\N$-grading. 
We note that if $\phi:P\to Q$ is surjective and 
$P,Q\in {\rm mod}_{\N}(V)$, then $\wt(Q)\subseteq \wt(P)$. \\

Let us start the proof of Main Theorem. \\
By the same arguments in the proof of Key Theorem, we know that 
$$ (\oplus_r U^{(r)})\boxtimes (\oplus_k W^{(k)})\cong 
\oplus_{r,k}\left( U^{(r)}\boxtimes W^{(k)}\right)$$
and so it is sufficient to prove Main Theorem for indecomposable $V$-modules $U$ and $W$. 

As we showed in the proof of Lemma 6, 
for $T=\oplus_{i=0}^{\infty}T_{(i)}\in {\rm mod}_{\N}(V)$, we have 
$C_m(T)\cap T_{(0)}=\{0\}$ by the definition of $C_m(T)$.  
We also note that if $A=\oplus_{i=0}^{\infty}A_{(i)}$ and 
$B=\oplus_{j=0}^{\infty}B_{(j)}$ are $\N$-graded $V$-modules and 
$\phi:A\to B$ is a surjective $V$-homomorphism, 
then $\phi(C_m(A))=C_m(B)$ since $v_{-m}\phi(a)=\phi(v_{-m}a)$ for 
$a\in A$ and $v\in V$. In particular, $\dim(B/C_m(B))\leq \dim(A/C_m(A))$. 

On the other hand, by Key Theorem, there is a number $f_m(U,W)\in \N$ such that 
$\dim T/C_m(T)<f_m(U,W)$ for any $(T,\CY)\in \CF(U,W)$. 
Therefore, there is $(S,\CJ)\in \CF(U,W)$ such that 
for any $(T,\CY)\in \CF(U,W)$ we have $\dim (T/C_m(T))\leq \dim (S/C_m(S))$. 
As we have shown, there is a finite set $\{a_i\in \C\mid i\in I\}$ such that 
$$\wt(S)\subseteq \cup_{i\in I}(a_i+\N). $$
We fix $(S,\CJ)$ and $\{a_i\mid i\in I\}$ for a while. \\

\noindent
{\bf Lemma 7}\quad {\it For any $(T,\CY)\in \CF(U,W)$, 
we have $\wt(T)\subseteq \cup_{i\in I} (a_i+\N).$}\\

\noindent
{\bf Proof}\quad 
For $(T,\CY)\in \CF(U,W)$, there is 
$(P,\CI)\in \CF(U,W)$ such that $(P,\CI)>(T,\CY)$ and $(P,\CI)>(S,\CJ)$. 
Since $(P,\CI)>(T,\CY)$, we have $\wt(T)\subseteq \wt(P)$ and so 
it is sufficient to prove Lemma 7 for $(P,\CI)$. 
Therefore, we may assume $(T,\CY)>(S,\CJ)$. 
Let $\phi:T\to S$ be a surjection. In this case, 
since $\phi(C_m(T)) \subseteq C_m(S)$ and $\dim (T/C_m(T))\leq \dim (S/C_m(S))$, 
we have 
$\Ker(\phi) \subseteq C_m(T)$ and $\Ker(\phi)\cap T_{(0)}=\{0\}$. 
Therefore, we have 
$$\wt(T)\subseteq \wt(T_0)+\N \subseteq \wt(S)+\N \subseteq \cup_{i\in I}(a_i+\N) $$
as we desired. \\

We come back to the proof of Main Theorem. 
Since 
$$(U\boxtimes W)_r\subseteq \prod_{(T,\CY)\in \CF(U,W)} T_r$$ 
by the construction (3.1), 
$$\wt(U\boxtimes W)\subseteq \cup_{i\in I}(a_i+\N).$$ 
Namely, the weights of elements in $U\boxtimes W$ is bounded below and so 
$v_nw=0$ for a sufficiently large $n$ for $w\in U\boxtimes W$ and $v\in V$. 
Therefore, $U\boxtimes W$ is a (weak) $\N$-gradable $V$-module. 
It is also $C_1$-cofinite as a $V$-module by Key Theorem. 

The remaining thing is to show that $\CY^{U\boxtimes W}$ is a (logarithmic) 
intertwining operator. By the construction of $\CY^{U\boxtimes W}$, 
it has a form: 
$$\CY^{U\boxtimes W}(u,z)w=\sum_{i=0}^{\infty}\sum_{n\in \C} u_{(i,n)}z^{-n-1}\log^iz. \eqno{(4.1)}$$ 
We have to prove that powers of $\log z$ in (4.1) are bounded.  
Since the homogeneous subspaces of $U, W$ and $U\boxtimes W$ 
are of finite dimension by Lemma 6, 
$L(0)^{{\rm nil}}=L(0)-\wt$ acts nilpotently on every homogeneous 
spaces $U_{(n)},W_{(n)}$ and $(U\boxtimes W)_{(n)}$. 
Furthermore, 
since $L(0)^{{\rm nil}}$ commutes with all actions $v_k$ for $v\in V$ and $k\in \Z$ and  
$U\boxtimes W$ is $C_1$-cofinite as a $V$-module, there is an integer $N$ such that 
$(L(0)^{{\rm nil}})^N(U\boxtimes W)=0$. 
Similarly, we may assume $(L(0)^{{\rm nil}})^NW=0$ and $(L(0)^{{\rm nil}})^NU=0$ by 
taking $N$ large enough. 
From the $L(-1)$-derivative property for $\CY^{U\boxtimes W}$, we have 
$$ (i+1)u_{(i+1,n)}w=-L(0)^{{\rm nil}}(u_{(i,n)}w)
+(L(0)^{{\rm nil}}u)_{(i,n)}w+u_{(i,n)}(L(0)^{{\rm nil}}w),\eqno{(4.2)}$$
Therefore, $u_{(k,n)}w=0$ for $k\geq 3N$, $u\in U$ and $w\in W$, which 
implies that $\CY^{U\boxtimes W}$ is a (logarithmic) intertwining operator.  

This completes the proof of Main Theorem.
\prend

\section{Discussion}

In this section, we would like to study the set $\CO^{C1}$ of all $\N$-gradable 
$C_1$-cofinite modules. 
By the above theorems, if $U,W\in \CO^{C1}$, then $U\boxtimes W\in \CO^{C1}$. 
Furthermore, there is a surjective (logarithmic) intertwining operator 
$\CY^{U\boxtimes W}\in I\binom{U\boxtimes W}{U\, \, W}$. 
We choose and fix such an intertwining operator $\CY^{U\boxtimes W}$ 
for each pair $U,W\in \CO^{C1}$. 

As Huang has shown in \cite{H}, the associativity of 
products of intertwining operators 
among $C_1$-cofinite $\N$-gradable $V$-modules are 
given by the expansion of solutions of the same differential 
equations in the suitable regions. 
Namely, if $U, W, T\in {\rm mod}_{\N}(V)$ are $C_1$-cofinite as $V$-modules,  
then by expanding solutions of differential equations on 
$|z_1|>|z_2|>0$, we know that there is an intertwining operator 
$\displaystyle{\CY\in I\binom{(U\boxtimes W)\boxtimes T}{U, \quad W\boxtimes T}}$ such 
that 
$$\langle p', \CY^{(U\boxtimes W)\boxtimes T}
(\CY^{U\boxtimes W}(u,z_1-z_2)w, z_2)t\rangle 
=\langle p', \CY(u,z_1)\CY^{W\boxtimes T}(w,z_2)t\rangle \eqno{(5.1)}$$
for $p'\in ((U\boxtimes W)\boxtimes T)', u\in U, w\in W$ and $t\in T$. 
We note that the right side in (5.1) was usually expressed by a linear sum, but 
we can express it with only one $\CY$ because of the universality property of $\CY^{W\boxtimes T}$.  
By the universality of $U\boxtimes (W\boxtimes T)$, 
there is a unique homomorphism 
$$\phi: U\boxtimes (W\boxtimes T)\to 
(U\boxtimes W)\boxtimes T$$ 
such that $\phi\circ\CY^{U\boxtimes (W\boxtimes T)}=\CY$. 
We also apply this argument to the opposite site. \vspace{2mm}\\
Namely, there is a homomorphism 
$\varphi: (U\boxtimes W)\boxtimes T\to U\boxtimes (W\boxtimes T)$ 
such that 
$$\langle p", (\varphi\circ\CY^{U\boxtimes W)\boxtimes T}) (\CY^{U\boxtimes W}(u,z_1-z_2)w, z_2)t\rangle 
=\langle p", \CY^{U\boxtimes (W\boxtimes T)}(u,z_1)\CY^{W\boxtimes T}(w,z_2)t\rangle$$
for $p"\in (U\boxtimes(W\boxtimes T))'$. 
In this case, since \vspace{1mm}
$$\begin{array}{l}
\displaystyle{\langle p', 
(\phi\varphi)\left(\CY^{(U\boxtimes W)\boxtimes T}\left(\CY^{U\boxtimes W}(u,z_1-z_2)w, z_2\right)t \right)\rangle}\vspace{1mm}\cr
\mbox{}\qquad=\displaystyle{\langle p', 
\phi\left(\varphi\left(\CY^{(U\boxtimes W)\boxtimes T}(\CY^{U\boxtimes W}(u,z_1-z_2)w, z_2)t \right)\right)\rangle}\vspace{1mm}\cr
\mbox{}\qquad=\displaystyle{\langle p', 
\phi(\varphi\circ\CY^{(U\boxtimes W)\boxtimes T}(\CY^{U\boxtimes W}(u,z_1-z_2)w, z_2)t )\rangle}\vspace{1mm}\cr
\mbox{}\qquad
=\displaystyle{\langle p',\phi(\CY^{U\boxtimes(W\boxtimes T)})(u,z_1)\CY^{W\boxtimes T}(w,z_2)t)\rangle }\vspace{1mm}\cr
\mbox{}\qquad
=\displaystyle{\langle p',\phi\circ\CY^{U\boxtimes(W\boxtimes T)})(u,z_1)\CY^{W\boxtimes T}(w,z_2)t\rangle }\vspace{1mm}\cr
\mbox{}\qquad =\displaystyle{\langle p', \CY^{(U\boxtimes W)\boxtimes T}(\CY^{U\boxtimes W}(u,z_1-z_2)w, z_2)t\rangle }\vspace{2mm}
\end{array}$$
for any $p'\in ((U\boxtimes W)\boxtimes T)'$, $u\in U$, $w\in W$ and $t\in T$. 
Therefore, $\phi\circ \varphi=1$. Similarly, we have $\varphi\circ\phi=1$ and so $\phi$ is an isomorphism. 
Namely, we have proved the following:\\

\noindent
{\bf Proposition 8}\quad {\it 
For $U,W,T\in \CO$, there is the unique isomorphism 
$$\phi:U\boxtimes (W\boxtimes T)\to (U\boxtimes W)\boxtimes T$$ 
satisfying
$$\langle p', \CY^{(U\boxtimes W)\boxtimes T}
(\CY^{U\boxtimes W}(u,z_1-z_2)w, z_2)t\rangle 
=\langle p', \phi\left(\CY^{U\boxtimes(W\boxtimes T)}(u,z_1)\CY^{W\boxtimes T}(w,z_2)t\right)\rangle $$
for $p'\in ((U\boxtimes W)\boxtimes T)', u\in U, w\in W$ and $t\in T$. }\\

At last, we would like to note one more thing. 
Although we have treated all logarithmic intertwining operators, 
if we restrict intertwining operators into formal $\C$-power series, that is, 
if we consider 
$$\begin{array}{l}
 \CF_{fp}(U,W)=\left\{(T,\CY)\in \CF(U,W)\mid \mbox{\begin{tabular}{l}
$\CY(u,z)w$ is a formal $\C$-power series \\
for any $u\in U$ and $w\in W$ \end{tabular}} \right\} \cr
\end{array} $$
then $\CF_{fp}(U,W)$ is still a direct set and so 
we are able to consider its projective limit, that is, 
another fusion product $U\boxtimes_{fp}W$ in this sense.

\end{document}